\documentclass{amsart}
\usepackage{array}
\usepackage{tikz}
\usepackage{xcolor}
\usepackage{moreverb}
\usepackage{amssymb,amsmath,amsthm}
\usepackage{mathtools}
\usepackage{enumitem}
\usepackage{tikz}
\usepackage{tikz-cd}
\usepackage{graphicx}
\usepackage{ulem}
\usepackage{url}
\usepackage{hyperref}

\usepackage[english]{babel}
\usepackage[latin1]{inputenc}
\usepackage{times}
\usepackage[T1]{fontenc}

\definecolor{olive}{rgb}{0.3, 0.4, .1}
\definecolor{fruitsalad}{rgb}{0.33, 0.67, 0.33}

\newcommand{\CC}{\mathbb{C}} 
\newcommand{\ZZ}{\mathbb{Z}} 
\newcommand{\NN}{\mathbb{N}} 
 
\newcommand{\Nz}{\mathbb{N}_0} 

\newcommand{\torus}{\mathbb{T}}
\newcommand{\disp}{\displaystyle}
\DeclarePairedDelimiter{\bracepair}{\lbrace}{\rbrace}
\DeclarePairedDelimiter{\anglepair}{\langle}{\rangle}
\DeclarePairedDelimiter{\vertpair}{\vert}{\vert}
\DeclarePairedDelimiter{\Vertpair}{\Vert}{\Vert}

\newcommand{\abs}[1]{\vertpair*{#1}}

\newcommand{\card}[1]{\vertpair{#1}}

\newcommand{\lpdualpair}[2]{\anglepair{#1, #2}}
\newcommand{\lpnormsh}[2]{\Vertpair{#1}_{#2}}

\newcommand{\spec}{\sigma}

\renewcommand{\Im}{\operatorname{Im}}

\newcommand{\unitdisk}{\mathcal{D}}
\newcommand{\clos}[1]{\overline{#1}}

\newtheorem{thm}{Theorem}

\newtheorem{clm}[thm]{Claim}

{ \theoremstyle{remark} }
{ \theoremstyle{definition} }
{ \theoremstyle{remark} \newtheorem{example}[thm]{Example}}

\title{Systems of Dilated Functions: completeness, minimality, basisness}
\author{Boris~Mityagin}
\address{The Ohio State University, Columbus, Ohio, USA}
\email{\url{mityagin.1@osu.edu}\\
\url{boris.mityagin@gmail.com}}

\begin{document}
\begin{abstract}
We discuss completeness, minimality, and basisness, in $L^2[0, \pi]$ and $L^p[0, \pi]$, $p \neq 2$, of dilated systems $u_n(x) = S(nx)$, $n \in \NN$, where $S$ is a trigonometric polynomial 
\[
S(x) = \sum_{k = 0}^m a_k \sin(kx), \quad a_0 a_m \neq 0.
\] 

We will present some results and mention a few unsolved questions. \nocite{*}
\end{abstract}
\maketitle

\normalem
\everymath{\displaystyle}

We announce a series of results which complement and extend the research of \cite{GNN} -- \cite{MitZero}.  The proofs will be given elsewhere.  

Define the polynomial
\[
a(z) = \sum_{j = 0}^m a_j z^j,
\]
and sets
\[
\begin{array}{c c c c c c c}
Z(a) & = & F^{-} & \cup & F^0 & \cup & F^+\\
a(\alpha) = 0 & & \abs{\alpha} < 1 && \abs{\alpha} = 1 && \abs{\alpha} > 1.
\end{array}
\]
The isometry $T: f(x) \to f(2x)$ has spectrum $\spec(T) = \clos{\unitdisk}, \quad \unitdisk = \bracepair{\zeta \in \CC: \abs{\zeta} < 1}$.  Then $u_n = a(T) \bracepair{\sin(nx)}$.  We factorize $a(z) = a^{-}(z) a^0(z) a^+(z)$.  We note that $a^+(T)$ is invertible, for letting $B = (a^{+}(T))^{-1}$ ,
\[
\begin{split}
B & = \prod_{\alpha \in F^+} (\alpha - T)^{-1}  = \frac{1}{2\pi i} \int\limits_{\abs{z} = 1 + \delta} R(z, t) \frac{dz}{a^+(z)} \, ,\\
1 + 2 \delta & = \min \bracepair{\abs{\alpha} : \alpha \in F^+}.
\end{split}
\]
Let $v_n = B u_n = a^{-}(T) a^{0}(T) \bracepair{\sin(nx)}$.

\begin{clm} \mbox{}\\
\begin{enumerate}[label = (\alph*)]
\item The system $U = \bracepair{u_n}$ is a basis in $L^2[0, \pi]$ if and only if $F^- \cup F^0 = \emptyset$.
\item If $F^- \neq \emptyset$ the system is not complete.
\item For any $a$ the system $U$ is minimal.  
\end{enumerate} 
\end{clm}

So assume $(\star)$ $F = F^0 \neq \emptyset$, i.e., all roots $\alpha$, $a(\alpha) = 0$, are in $\mathbb{T} = \bracepair{z \in \CC: \abs{z} = 1}$.  
\begin{clm}
Under $(\star)$
\[
a(z) = a_m \prod_{\alpha \in F^0} (\alpha - z)^{\mu(\alpha)}.  
\]
Put $\kappa^* = \max \bracepair{\mu(\alpha) - 1: \alpha \in F^0}$.  The system $U$ is complete, and minimal, i.e., $\exists \bracepair{\Phi_k}$, $\lpdualpair{\Phi_k}{u_n} = \delta_{kn}$, and 
\[
\lpnormsh{\Phi_k}{q} \asymp \left( \log k \right)^{\kappa^* + 1/2},
\] 
so $U$ is not a basis in $L^2$ (or $L^p$, $1 < p < \infty$).  
\end{clm}
Now we go to the ``multi-frequency'' case.  We again set $U = \bracepair{u_n(x)}$ and $u_n(x) = S(nx)$, $0 \leq x \leq 2 \pi$, 
\[
\begin{split}
S(x) &= \sum_{j \in J} a_j \exp (i j x),\, \, \card{J} < \infty\\
& = \sum_{\substack{ \alpha \in \Nz^m \\ \alpha \in K, \, \, K \subseteq \Nz^m }} a(\alpha) \exp \left( i \left[ \prod_{j = 1}^m p_j^{\alpha_j} \right] x \right), \, \,  \card{K} < \infty,
\end{split}
\]
where $\bracepair{p_j}_{j = 1}^m$ is a set of primes.  Put
\[
A(\omega) = \sum_{\alpha \in K} a(\alpha) w^{\alpha}, \quad w^{\alpha} = \prod_{j = 1}^m w_j^{\alpha_j}.
\] 
and 
\[ 
\begin{split}
\NN(\omega) &= \bracepair{\omega \cdot p^{\alpha} : \alpha \in \Nz^m}\\
\omega \in \Omega & = \bracepair{q \in \NN: q \text{ does not have factors } p_j, \, \, 1 \leq j \leq m}.
\end{split}
\]
Let 
\[
T_j: f(x) \mapsto f(p_j x), \quad 1 \leq j \leq m.
\]
Consider
\[
E \cong \ell^2 \simeq H^2(\unitdisk) \simeq \ell^2 \left( \Omega; H^2(\unitdisk^m) \right).
\]
All $E(\omega) = \Im Q(\omega)$, $\omega \in \Omega$, are invariant with respect to the isometries $T_j$, $1 \leq j \leq m$, i.e. multiplication by $w_j$ in $H^2(\unitdisk^m)$. 
Certainly,
\[
\Vertpair{f}^2 = \sum_{\omega \in \Omega} \Vertpair{Q(\omega) f}^2.
\]

Now all the questions about $U$ become the questions about the system
\[
 V = \bracepair{v(\alpha)}_{\alpha \in \Nz^m} ; \quad v(\alpha)(w) = A(w) w^{\alpha} \quad \text{ in } H^2(\unitdisk^m).
\]
\begin{clm}
If 
\begin{equation} \label{eq:zoutsideOK} \tag{$\star \star$}
Z(A) \cap \clos{\unitdisk^m} = \emptyset,
\end{equation}
then $A(T)^{-1} = B$ is well-defined and $\bracepair{v(\alpha)}$ is a (Riesz) basis and $U$ is a (Riesz) basis as well.

If $V$ (or $U$) is a Riesz basis then \eqref{eq:zoutsideOK} holds.  
\end{clm}

\begin{clm}
The system $U$ is minimal if $a(0) \neq 0$.
\end{clm}

Indeed, with 
\[
\lpdualpair{f}{g} = \frac{1}{(2\pi)^m} \int\limits_{\torus^m} f(w) g(w) \, \mathrm{d}^mt, \quad w_j = e^{it_j}, \quad 1 \leq j \leq m,
\]
(no bar, no conjugation), $\lpdualpair{w^{-\tau}}{w^{\alpha}} = \delta(\alpha, \tau), \quad \forall \alpha, \tau \in \ZZ^m$.  Then
\[
\frac{1}{A(w)} = \sum_{\sigma \in \Nz^m} b(\sigma) w^{\sigma}, \quad \frac{w^{-\tau}}{A(w)} = \sum_{\sigma \in \Nz^m} b(\sigma) w^{\sigma - \tau},
\]
but if $\sigma - \tau \leq 0$ does not hold then
\[
\lpdualpair{w^{\sigma - \tau}}{w^{\alpha}} = 0, \quad \forall \alpha \in \Nz^m ,
\]
so put $\Phi_t(w) = \sum_{\sigma \leq \tau} b(\sigma) w^{\sigma - \tau}$. 
This \emph{finite} sum is well-defined, and 
\[
\lpdualpair{\Phi_{\tau}}{v(\alpha)} = \delta(\alpha, \tau), \quad \text{ for all } \alpha, \tau \in \Nz^m.
\]

\[
Z(A) \cap \clos{\unitdisk^m} = Z(A) \cap \torus^m.
\]

Completeness implies uniqueness of the system $\Phi_{\tau}$ and the fact that for 1D projections $P_{\tau} = \lpdualpair{\bullet}{\Phi_{\tau}} v(\tau)$,
\[
\Vertpair{P_{\tau}} = \Vertpair{\Phi_{\tau}} \cdot \Vertpair{v(\tau)} \asymp \Vertpair{\Phi_{\tau}}.
\]
But
\[
\Vertpair{\Phi_{\tau}}^2 = \sum_{\sigma \leq \tau} \abs{b(\sigma)}^2, \quad B(w) = \frac{1}{A(w)},
\]
so these norms are uniformly bounded if and only if
\[
\frac{1}{A(w)} \in H^2(\unitdisk^m), \text{ or } \frac{1}{P(t)} \in L^1(\torus^m), \text{ where } P(t) = \abs{A(e^{it})}^2.
\]

\begin{clm}
If $m \leq 3$ and $Z(A) \cap \clos{\unitdisk^m} \neq \emptyset$, then
\[
\frac{1}{A(w)} \not\in H^2(\unitdisk^m).
\]
Under the same assumptions, $V$ or $U$ is \emph{NOT} a basis.  
\end{clm}

For $m \geq 4$, it could happen that $\disp \frac{1}{A(w)} \in H^2(\unitdisk^m)$.  
\begin{example}
Fix $c_k > 0$, $1 \leq k \leq m$, with $\sum_{k = 1}^m c_k = 1$, and define
\begin{equation} \label{eq:aobstructdef}  \tag{$E^*$}
A(w) = 1 - \sum_{k = 1}^m c_k w_k.
\end{equation}
\end{example}

Then 
\[
\begin{split}
P(t) &= \abs{\sum_{k = 1}^m c_k 2 \sin^2 \left( \frac{t_k}{2} \right)}^2 + \abs{\sum_{k = 1}^m c_k \sin(t_k)}^2\\
& \asymp r^4 + \abs{\ell(t)}^2, \quad \ell(t) = \sum_{k = 1}^m c_k t_k, \quad \vert t \vert \ll 1.  
\end{split} 
\]
We note that $\disp \int\limits_{\abs{\zeta} \leq \delta} \frac{d\zeta_0 \, d\zeta_1 \dotsc d\zeta_{m-1}}{\zeta_0^2 + \zeta_1^4 + \dotsb + \zeta_{m-1}^4} < \infty$ if and only if $m \geq 4$, since it is within a constant multiple of
\[
\int_0^{\delta} \int_0^{\rho^2} \frac{d\zeta \,  \rho^{m-2} \, d\rho}{\zeta^2 + \rho^4} = \int_0^{\delta} \int_0^1 \frac{d\eta \,  \rho^{m-4} \, d\rho}{1 + \eta^2}.
\]

Recall that $v(\alpha) = A(w) w^{\alpha}$, $\alpha \in \Nz^m$.  Instead of asking whether this system is a basis in $H^2(\unitdisk^m)$, or $\ell^2(\Nz^m)$ we can move to weighted $H^2(\unitdisk^m; P)$, or $L^2(\torus^m; P)$ and ask whether $\bracepair{w^{\alpha}}$ is a basis.  

\begin{clm}
In the case \eqref{eq:aobstructdef}, $m \geq 4$, for the partial sums $\Sigma(\tau) f = \sum_{\alpha \leq \tau} \lpdualpair{\Phi_{\alpha}}{f} v_{\alpha}$, the norms $\Vertpair{\Sigma(\tau)}$ are not bounded.  
\end{clm}
The proof is based on the multi-dimensional $A_2$ Muckenhoupt condition (\cite{KazLiz}, \cite{Moen}).

But if we go to the original system
\[
v_{\alpha}(x) = \sum_{k = 0}^K a_k \exp \left( i [p^{\alpha}] kx \right), \quad p = (p_j)_{j = 1}^m, \quad \alpha \in \Nz^m,
\]
its linear ordering would fit to a monotone arrangement of the multi-index sequence $\bracepair{p^{\alpha}}$, or its linear ordering by monotonicity of the linear form $M(\alpha) = \sum_{j = 1}^m \alpha_j \log p_j$.
It leads us to the question on the boundedness of the projection $Q_M$.

\[
\begin{split}
\exp i (\alpha, t) & \to \text{ the same, if } M(\alpha) \geq 0\\
& \to 0 \text{ if} M(\alpha) < 0.  
\end{split}
\] 
in the weighted $L^2(\torus^m; P)$, $m \geq 4$, say, 
\[
P(t) = \left( \sum_{j = 1}^m t_j \right)^2 + \left( \sum_{j = 1}^m t_j^2 \right)^2.
\]
$M(\alpha)$ is ``essentially irrational,'' i.e., the coefficients $\mu_j^* = \log p_j, \quad 1 \leq j \leq m$, of the linear function $M(y) = \sum_{j = 1}^m \mu_j y_j$ are rationally independent.  

If all $\mu_j$ were rational $Q_M$ would be equivalent to the case $M_0(y) = y_1$; then known $A_2$ conditions are applicable and $Q_{M_0}$ and $Q_M$ are bounded, for any $m$.  

If $Q_{M^*}$ were bounded, $m \geq 4$ (I do not believe so) we would have Babenko--type Shauder (but not Riesz) basis in $H^2(\unitdisk^m)$, or even in $L^2(\torus^m)$.

The announced results were presented at the Workshop on Operator Theory, Complex Analysis, and Applications (WOTCA), Coimbra, Portugal, June 21 -- 24, 2016.  

\bibliographystyle{unsrt}
\bibliography{DilationsShort}

\end{document}